\documentclass[11pt]{article}
\usepackage{amsmath,amsfonts,amssymb,amsthm,xypic,epsfig,setspace}
\usepackage[pdftex]{hyperref}

\theoremstyle{plain}
\newtheorem{theorem}{Theorem}[section]
\newtheorem{proposition}[theorem]{Proposition}

\theoremstyle{definition}
\newtheorem{lemma}[theorem]{Lemma}

\newtheorem{cor}[theorem]{Corollary}
\newtheorem{example}[theorem]{Example}

\newenvironment{renumerate}%
{%
\begin{enumerate}}%
{\end{enumerate}%
}%

\newenvironment{definition}%
{\vskip6pt%
\noindent%
{\bf Definition.}}%
{\vskip6pt}

\newenvironment{remark}%
{\vskip6pt%
\noindent%
{\bf Remark: }}%
{\vskip6pt}

{\vskip6pt%
\noindent%
{\it Remarks}. %
\begin{renumerate}}%
{\end{renumerate}\vskip6pt}

\newcommand{\R}{\text{${\mathbb R}$}}
\newcommand{\C}{\text{$\mathbb C$}}

\newcommand{\Z}{\text{$\mathbb Z$}}

\newcommand{\M}{\mathcal{M}}
\newcommand{\Hh}{\mathcal{H}}
\newcommand{\Kperp}{K^{\perp}}
\renewcommand{\frak}[1]{\text{$\mathfrak{#1}$}}
\newcommand{\J}{\mathcal{J}}
\newcommand{\G}{\mathcal{G}}
\newcommand{\E}{\mathcal{E}}

\newcommand{\End}{\mathrm{End}}
\newcommand{\Gl}{\mathrm{Gl}}

\newcommand{\lra}{\longrightarrow}

\newcommand{\gO}{\text{$\Omega$}}
\newcommand{\gf}{\text{$\varphi$}}

\newcommand{\Id}{\mathrm{Id}}

\newcommand{\del}{\text{$\partial$}}

\newcommand{\delbar}{\text{$\overline{\partial}$}}
\newcommand{\tensor}{\otimes}

\newcommand{\mc}[1]{\text{$\mathcal{#1}$}}
\newcommand{\mr}[1]{\text{$\mathrm{#1}$}}

\newcommand{\into}{\rightarrow}

\newcommand{\noqed}{\let\qed\relax}

\newcommand{\Gg}{\mathfrak{g}}

\newcommand{\gcs}{generalized complex structure}

\newcommand{\gcss}{generalized complex structures}

\newcommand{\gk}{generalized K\"ahler}
\newcommand{\gks}{generalized K\"ahler structure}

\newcommand{\gcy}{generalized Calabi--Yau}

\newcommand{\wrt}{with respect to}
\renewcommand{\iff}{if and only if}
\newcommand{\Ann}{\mathrm{Ann}}
\newcommand{\IP}[1]{\langle #1 \rangle}

\numberwithin{equation}{section}
\begin{document}
\title{Generalized complex geometry and T-duality}
\author{
Gil R. Cavalcanti\thanks{\tt g.r.cavalcanti@uu.nl}~ and Marco Gualtieri \thanks{\tt mgualt@math.toronto.edu}
}
\date{}

\maketitle

\abstract{ We describe how generalized complex geometry, which interpolates between complex and symplectic geometry, is compatible with T-duality, a relation between quantum field theories discovered by physicists.  T-duality relates topologically distinct torus bundles, and prescribes a method for transporting geometrical structures between them.  We describe how this relation may be understood as a Courant algebroid isomorphism between the spaces in question. This then allows us to transport Dirac structures, generalized
Riemannian metrics, generalized complex and generalized K\"ahler structures, extending the
\emph{Buscher rules} well-known to physicists.  Finally, we re-interpret T-duality as a Courant reduction, and explain that T-duality between generalized complex manifolds may be viewed as a generalized complex submanifold (D-brane) of the product, in a way that
establishes a direct analogy with the Fourier-Mukai transform.
}


\section*{Introduction}
T-duality is an equivalence between quantum field theories with very
different classical descriptions; for example type IIA and IIB
string theory are T-dual when compactified on a circle. The precise
relationship between T-dual Riemannian structures was first
understood by Buscher in \cite{Bu87} and was developed further by
Ro\v{c}ek and Verlinde in \cite{RoVe92}.  It was realized that in
order to phrase T-duality geometrically, one had to consider the
interplay between the Neveu-Schwarz 3-form flux $H$, a closed 3-form
with integral periods which entered the sigma model as the
Wess-Zumino term, and the topology of the sigma model target. The
precise relation between this 3-form flux and the topology of the
T-dual spaces has recently been given a clear description by
Bouwknegt, Evslin and Mathai in~\cite{BEM03} and it is their
topological approach which we shall use as a basis to study the
geometry of T-duality.

In this paper we explore and expand upon the realization
in~\cite{Gu03} that T-duality transformations can be understood in
the framework of \emph{generalized geometrical structures}
introduced by Hitchin in~\cite{Hi03}.  In this formalism, one
studies the geometry of the direct sum of the tangent and cotangent
bundles of a manifold.  This bundle is equipped with a natural
orthogonal structure as well as the \emph{Courant bracket}, an
analog of the Lie bracket of vector fields, which depends upon the
choice of a closed 3-form.  In particular, an integrable orthogonal
complex structure on this bundle, or generalized complex structure,
is an object which encompasses complex and symplectic geometry as
extremal special cases.  As we shall see, T-duality can be viewed as
an isomorphism between the underlying orthogonal and Courant
structures of two possibly topologically distinct manifolds. It can
therefore be used to transport a generalized complex structure from
one manifold to the other, and in so doing, complex and symplectic
structures on the two manifolds may be interchanged. This helps us
to more fully understand the proposal of~\cite{SYZ96} that mirror
symmetry between complex and symplectic structures on Calabi-Yau
manifolds can be understood as an application of T-duality.

The action of T-duality on generalized complex structures was
implicitly observed in~\cite{FMT03}, where both complex and
symplectic structures in 6 dimensions were interpreted as spinors
for $CL(6,6)$, a natural consideration from the point of view of
supergravity. However, without the formalism of generalized complex
structures, the intermediate geometrical structures were not
recognized. Once the connection with generalized geometry was
understood, several works appeared
\cite{GMPT04,GLMW03,CGJ06,LMTZ04,LRUZ04,ZU05,Zu04} which provide a
physical motivation and justification for the use of generalized
complex structures to understand mirror symmetry. From a
mathematical point of view, Ben-Bassat~\cite{BeB06,BeB06b} explored the
action of T-duality on generalized complex structures on vector
spaces and flat torus bundles, where one does not consider the
3-form flux $H$ and therefore restricts the topological type of the
bundles in question.

While we treat the most general case of T-duality of circle bundles
with 3-form flux, it is important to clarify that for higher rank
affine torus bundles, we only consider 3-forms $H$ for which $i_X
i_Y H=0$ for $X,Y$ tangent to the fibres. Mathai and
Rosenberg~\cite{MaRo04} have shown that without this restriction,
the T-dual manifold may be viewed as a noncommutative space.  While
this may also have an interesting interpretation in terms of
generalized geometry, we do not explore it here.

{\bf Acknowledgements:} We wish to thank Nigel Hitchin for many
helpful discussions and insights.  We are grateful to Christopher
Douglas and Lisa Jeffrey for assistance.  This research was supported
in part by EPSRC, NSERC, NWO  as well as the Fields Institute.  We would also like to thank Robert Kotiuga for creating and organizing, with support from CRM, NSF, and CMI, a wonderful and memorable conference in honour of Raoul Bott, at which some of this work was presented.  

\section{Generalized complex geometry}

Given a real $n$-dimensional manifold $M$ and a closed 3-form $H\in \Omega^3(M,\R)$, the sum of the tangent and cotangent bundles of $M$, which we denote by $T\oplus T^*$, is endowed with a natural symmetric bilinear form of signature  $(n,n)$
$$\IP{X + \xi,Y + \eta} = \tfrac{1}{2}(\eta(X)  +\xi(Y)),$$
as well as a bracket operation, the {\it Courant bracket}, given by 
$$[X+\xi,Y +\eta] = [X,Y]+ \mc{L}_X \eta -i_Yd\xi + i_Yi_X H.$$ 
Elements of  $T \oplus T^*$ act naturally on the space of forms by the parity reversing map
$$ (X +\xi) \cdot \rho  = i_X\rho + \xi\wedge \rho.$$
This action extends to an action of the Clifford algebra bundle $\mr{Cl}(T\oplus T^*)$ since,  for $v \in T \oplus T^*$, we have
$$ v\cdot(v\cdot \rho) = \IP{v,v}\rho,$$
rendering  $\wedge^{\bullet}T^*M$ into an irreducible Clifford module for $\mr{Cl}(T\oplus T^*)$.  Essentially, we may view differential forms on $M$ as spinors for the metric bundle $T\oplus T^*$. 

The closed 3-form $H$ defines a twisted de Rham differential $d_H$, defined by $d_H\rho = d\rho + H\wedge \rho$, which is related to the Courant bracket in the same way that the usual de Rham differential is related to the Lie bracket, i.e.  
\begin{equation}\label{eq:twisted bracket}
[v_1,v_2]_{H}\cdot \rho =  [[d_H,v_1],  v_2 ]\cdot \rho.
\end{equation}

An important symmetry of $T \oplus T^*$ is given by the action of  2-forms by {\it $B$-field transforms}. Given  $B \in \Omega^2(M)$, we view it as a map $B:T \into T^*$ which exponentiates to  
$$e^{-B}( X + \xi) = X + \xi - B(X).$$
This orthogonal map  relates the $H$-Courant bracket to the $H-dB$-Courant bracket:
$$[X+ \xi -B(X),Y+\eta-B(Y)]_H = [X+\xi,Y+\eta]_{H-dB} - B([X,Y]).$$ 
In the particular case in which $B$ is a closed form, it defines an orthogonal Courant automorphism of $T\oplus T^*$.

Using the action of the Clifford algebra, the map on forms corresponding to $B$-field transforms is given by $\rho \mapsto e^B \wedge \rho$, as we have the identity
$$(X + \xi - B(X))\cdot e^B\wedge \rho = e^B\wedge((X+\xi)\cdot \rho).$$
Therefore, a 2-form defines both an orthogonal map of $T\oplus T^*$ and an automorphism of the Clifford module $\wedge^{\bullet}T^*M$, compatible with the Clifford action.

An additional structure which is present on the Clifford module is the \emph{Mukai pairing}.  It is a bilinear form on $\wedge^\bullet T^*M$ with values in the volume forms $\wedge^n T^*M$, given by 
$$(\rho_1, \rho_2) := (\sigma(\rho_1)\wedge \rho_2)_{top},$$
where $(\cdot)_{top}$ indicates taking the top degree component and $\sigma$ is the linear map which acts on decomposable forms by
$$ \sigma(\xi_1 \wedge \cdots \wedge \xi_k) = \xi_k \wedge \cdots \wedge \xi_1, \qquad \xi_i \in T^*.$$
The Mukai pairing is compatible with the Clifford action, in the sense
$$(v\cdot \rho_1,v\cdot \rho_2) = \IP{v,v} (\rho_1,\rho_2)\qquad  \mbox{for all } v \in T \oplus T^*.$$
As a consequence, it is $\mathrm{Spin}_0$-invariant, in particular, for $B \in \Omega^2(M)$ we have
$$(e^B\wedge \rho_1,e^B\wedge \rho_2) = (\rho_1,\rho_2).$$

\subsection{Generalized geometrical structures}
Let $(M,H)$ be an manifold equipped with a closed real 3-form $H$.  We now review a series of geometrical structures on $T\oplus T^*$ which, as we explain in the subsequent section, may be transported across a T-duality relation between manifolds. 

\begin{definition}
A {\it Dirac structure} on $(M,H)$ is a subbundle $L \subset T \oplus T^*$  whose fibers are maximal isotropic subspaces and such that $\Gamma(L)$ is closed under the $H$-Courant bracket.
\end{definition}

\begin{definition}
A {\it \gcs}  on $(M^{2n},H)$ is a complex structure  $\J$ on $T \oplus T^*$, orthogonal with respect to the natural pairing,  and whose $+i$-eigenbundle, $L \subset (T\oplus T^*)\tensor \C$, is closed under the $H$-Courant bracket.
\end{definition}

Compatibility with the natural pairing implies that the $+i$-eigenbundle $L$ of a \gcs\  is a maximal isotropic subspace of $(T \oplus T^*)\otimes\C$. Since $L$ determines $\J$ completely, we see that a \gcs\ is a complex Dirac structure $L$ satisfying $L \cap \overline{L} = \{0\}$.

A \gcs\ may be described in terms of the Clifford module $\wedge^\bullet T^*M$ by specifying the unique complex line subbundle  $K\subset (\wedge^{\bullet}T^*M)\otimes\C$ such that $L$ is its Clifford annihilator, i.e.   
$$ L = \{v \in (T \oplus T^*)\otimes \C:v \cdot K =0\}.$$
 If $\rho$ is a nonvanishing local section of this line bundle, each of the conditions imposed on $L$ corresponds to a constraint on $\rho$, as follows \cite{Gu03}. The fact that $L$ is a maximal isotropic satisfying $L \cap \overline{L} = \{0\}$ translates to
$$\rho = e^{B+i\omega}\wedge \gO \qquad \mbox{and} \qquad  \Omega\wedge \overline{\Omega} \wedge (2i\omega)^{n-k} = (\rho,\overline{\rho}) \neq 0,$$
for $B,\omega$ real 2-forms, $\gO$ a decomposable complex $k$-form for some $k$, and $2n=\dim(M)$.  Also, $L$ is closed for the $H$-Courant bracket bracket \iff\
there is a local section $v\in\Gamma(T\oplus T^*)$ such that $$d_H\rho = v \cdot \rho.$$
\begin{definition}
The \emph{canonical bundle} of a generalized complex structure is the line subbundle  $K \subset (\wedge^{\bullet}T^*M)\otimes\C$ annihilated by $L = \ker (\J - i\Id)$.
\end{definition}

The {\it type} of a generalized complex structure at a point $p$ is the dimension of $\ker(\pi:L \into TM\otimes\C)$. If $\rho = e^{B+i \omega} \wedge \Omega$ generates the canonical bundle at $p$, then the type at $p$ is simply the degree of $\Omega$.

\begin{example}\label{ex:complex}
Any complex structure $I$ on a manifold  $M^{2n}$ with background 3-form $H=0$ gives rise to a \gcs\ $\J_I$. In the decomposition $T\oplus T^*$ we may write $\J_I$ as a block matrix:
$$ \J_I = \begin{pmatrix}
-I & 0 \\ 0 & I^*
\end{pmatrix}$$
The $+i$-eigenspace for this structure is $L = T^{0,1}M \oplus T^{*1,0}M$ and the canonical bundle is $\wedge^{n,0}T^*M$. This structure has type $n$ everywhere.
\end{example}

\begin{example}\label{ex:symplectic}
If $H=0$, a symplectic structure $\omega$ on $M$ also gives rise to a \gcs:
$$ \J_{\omega} = \begin{pmatrix}
0 & -\omega^{-1} \\ \omega & 0
\end{pmatrix}.$$
The $+i$-eigenspace of $\J_{\omega}$  is
$$L = \{X - i \omega(X): X \in TM\otimes\C\},$$
and the canonical bundle is generated by $e^{i\omega}$,  a nowhere
vanishing closed section. This structure has type zero everywhere. 
\end{example}

\begin{definition}
A {\it generalized Riemannian metric} on a manifold $M$ is an orthogonal, self-adjoint bundle map $\mc{G} \in \mbox{End}(T\oplus T^*)$ for which $\langle \mc{G}v, v \rangle$ is positive definite.
\end{definition}

A generalized metric $\mc{G}$, being self adjoint and orthogonal, must satisfy $\mc{G} = \mc{G}^* = \mc{G}^{-1}$.  Therefore $\mc{G}^2=\Id$, and $T \oplus T^*$ splits as an orthogonal sum of $\pm 1$-eigenspaces $C_{\pm} \subset T\oplus T^*$. As $\mc{G}$ is positive definite, the natural pairing is $\pm$-definite on $C_{\pm}$. Conversely, a generalized metric is determined by specifying a maximal positive-definite subbundle  $C_+ \subset T \oplus T^*$.

The subspace $C_+\subset T\oplus T^*$ may be expressed as the graph of a linear map $A:T \into T^*$, i.e.  an element in $\tensor^2 V^*$.  Using the splitting  $\tensor^2T^* = Sym^2 T^* \oplus \wedge^2 T^*$, we may write $C_+$ as the graph of $b+g$, where $g$ is a symmetric 2-form and $b$ is skew.  A similar analysis for $C_-$ shows:
$$ C_\pm = \{ X + b(X, \cdot)\pm g(X,\cdot)\ :\ X \in T\}.$$
In this way, a generalized Riemannian metric on $T\oplus T^*$  is  equivalent to a choice of a usual Riemannian metric $g$ and 2-form $b\in\Omega^2(M,\R)$.

\begin{definition}
A generalized complex structure $\J$ is compatible with a generalized Riemannian metric $\G$ if they commute.  Such a pair $(\G,\J)$ is called a {\it generalized Hermitian structure}.
\end{definition}

The compatibility of $(\G,\J)$ in a generalized Hermitian structure implies that  $\J_2 = \G \J$ satisfies $\J_2^2 = -\Id$. The structure $\J_2$  is the analogue of an almost complex structure taming a symplectic form or of a K\"ahler form for a Hermitian manifold; in general, it is not integrable. The integrablity of $\J_2$ is the analogue of the K\"ahler condition.

\begin{definition}
A {\it \gk\ structure} on $(M,H)$ is a pair of commuting \gcss\ $\J_1$ and $\J_2$ such that $\mc{G} =-\J_1 \J_2$ is a generalized Riemannian metric.
\end{definition}
A \gcs\ $\J$ on the real $2n$-manifold $M$ induces a decomposition of the bundle $\wedge^{\bullet}T^*M\otimes\C$ into subbundles $U^k$, similar to the $(p,q)$-decomposition induced by a complex structure. The top degree component, $U^n$, coincides with the canonical line bundle $K$, and we define
$$U^{n-k} = \wedge^{k}\overline{L}\cdot  U^n.$$
The bundle $U^k$ can also be described as the $ik$-eigenbundle of the Lie algebra action of $\J$ (see \cite{Gu03}). In the case of a usual complex structure, one obtains
$$U^k = \oplus_{p-q=k}\wedge^{p,q}T^*M\otimes\C,$$
whereas a symplectic structure has eigenbundles
$$U^{k} = e^{i\omega}e^{\frac{-\omega^{-1}}{2i}}\wedge^{n-k}T^*M\otimes\C,$$
where $\omega^{-1}$ denotes the bivector inverse to $\omega$, acting via interior product (see \cite{Ca05}).

For a \gk\ manifold, both $\J_1$ and $\J_2$ give rise to decompositions of forms.  Since $\J_1,\J_2$ commute, we obtain a bigrading of forms:
$$ U^{p,q} = U_{\J_1}^p\cap U_{\J_2}^q, \qquad  \Omega^{\bullet}(M,\C) =\oplus_{p,q} U^{p,q}.$$

For a usual K\"ahler manifold $(I,\omega)$, this decomposition is {\it not} the ordinary $(p,q)$ decomposition determined by the complex structure, but is obtained from it in the following way:
\begin{equation}\label{eq:upqkahler}
U_I^{p-q}\cap U_{\omega}^{n-p-q} =
e^{i\omega}e^{\frac{-\omega^{-1}}{2i}}\wedge^{p,q}T^*M\otimes\C.
\end{equation}

The decomposition of forms provides an analogue of the Dolbeault complex for any \gcs (see \cite{Gu03}, Theorem 4.23):  if $\mc{U}^k$ is the sheaf of local sections of $U^k$, then integrability implies the inclusion
$$d_H (\mc{U}^k) \subset \mc{U}^{k-1} \oplus \mc{U}^{k+1}.$$
This induces a splitting $d_H = \del + \delbar$, where:
$$\del: \mc{U}^k \into \mc{U}^{k+1}, \qquad \delbar:\mc{U}^k \into \mc{U}^{k-1}.$$
For complex manifolds, these coincide with the usual Dolbeault partial derivatives.  For the further decomposition induced by a generalized K\"ahler structure, see~\cite{GuK}.

The information in a  \gks\ $(\J_1,\J_2)$ may actually be rephrased in terms of a bi-Hermitian structure first discovered by Gates, Hull, and Ro\v{c}ek~\cite{GHR}, in their study of sigma models with $N=(2,2)$ supersymmetry.  A bi-Hermitian structure consists of a Riemannian metric $g$ with two compatible complex structures $I_+,I_-$ which need not commute;  in~$\cite{GHR}$, the following additional constraints were imposed:
\begin{equation}\label{eq:gk bihermitian}
d^c_- \omega_- + d^c_+\omega_+=0 \qquad \mbox{and}\qquad dd^c_-\omega_- =0,
\end{equation}
where $\omega_\pm = gI_\pm$ are the (not necessarily closed) K\"ahler forms for $I_\pm$, and $d^c_\pm = i(\delbar_\pm - \del_\pm)$. 
In fact, such a bi-Hermitian structure may be obtained from any generalized K\"ahler structure in the following way.   Since the metric $\G,\J_1$ commute, $\J_1(C_\pm) = C_\pm$, and hence $\J_1$ induces complex structures on $C_\pm$. By projection to $TM$, we obtain integrable complex structures $I_{\pm}$ on $TM$, compatible with the metric $g$. The bihermitian structure  $(g,I_+,I_-)$  thus obtained also satisfies~\eqref{eq:gk bihermitian}, and furthermore $[H]=[d^c_-\omega_-]$ in $H^3(M,\R)$. These results, shown in~\cite{Gu03}, may be reversed, yielding the following result.
\begin{theorem}[Gualtieri \cite{Gu03}]\label{theo:bihermitiangk}
A bihermitian structure $(g, I_{\pm})$ induces a \gks\  if and only if \eqref{eq:gk bihermitian} holds. 
\end{theorem}

The final structure we wish to recall is that of {\it strong K\"ahler with torsion} (SKT) geometry.  A Hermitian manifold $(M,I,g)$ is SKT if $dd^c\omega= 0$, where $\omega=gI$.  The previous theorem shows that a \gks\ consists of a pair of SKT structures \wrt\ the same metric, and with opposite torsion terms.  

An equivalent definition of SKT structure, following the argument used above for generalized K\"ahler geometry,  is as follows.
\begin{definition}
A SKT structure on $(M,H)$ is a generalized metric $\G$, together with an orthogonal  complex structure $\mc{I}:C_- \into C_-$ whose $+i$-eigenspace is closed under the Courant bracket. 
\end{definition}
From this point of view, we obtain a Hermitian structure $(g,I_-)$ by projection $\pi:C_- \into TM$,  for which $dd^c_-\omega_- = 0$ and $[d_-^c\omega_-] = [H].$

\section{T-duality with NS flux}\label{top}

In this section we review the definition of the T-duality relation by Bouwknegt, Evslin, Mathai and Hannabuss \cite{BEM03,BHM03}.  T-duality is a relation between principal torus bundles $M$ and $\tilde{M}$, both over the same base $B$, and each of which is equipped with closed 3-forms $H, \tilde H$, known to physicists as background NS-fluxes.  We review the main result of those papers which concerns us, which states that the twisted differential complexes of torus-invariant forms are isomorphic for T-dual manifolds, i.e. $(\Omega_T^{\bullet}(M),d_H) \cong (\Omega_T^{\bullet}(\tilde{M}),d_{\tilde{H}})$.

\begin{definition}\label{def:t-dual} Let $M$ and $\tilde{M}$ be two principal $T^k$ bundles over a common base manifold $B$ and let $H\in \Omega_{T^k}^3(M)$ and $\tilde{H}\in \Omega_{T^k}^3(\tilde{M})$ be invariant closed forms. Consider  $M \times_B \tilde{M}$, the fiber product of $M$ and $\tilde{M}$, so that we have the following diagram
\begin{equation}\label{eq:T-duality diagram}
\xymatrix{ & (M\times_{B}\tilde{M}, p^*H - \tilde{p}^*
\tilde{H})\ar[dl]_p\ar[dr]^{\tilde{p}} &
\\ (M,H)\ar[dr]_\pi& &(\tilde{M},\tilde{H})\ar[dl]^{\tilde{\pi}} \\ &B & \\}
\end{equation}
We say that $M$ and $\tilde{M}$ are  {\it T-dual} if 
\begin{equation}\label{eq:T-duality condition}
p^*H - \tilde{p}^*\tilde{H} =dF,
\end{equation}
 for some $T^{2k}$-invariant 2-form $F \in \Omega^2(M \times_B \tilde{M})$ such that
 \begin{equation}\label{eq:nondegenerate}
 F:\frak{t}_M^k \otimes \frak{t}_{\tilde{M}}^k \into \R \qquad\mbox{is nondegenerate,} 
 \end{equation}
where $\frak{t}^k_M$ is the tangent space to the torus fiber of $M\times_{B}\tilde{M} \into \tilde{M}$ and $\frak{t}_{\tilde{M}}^k$ the tangent to the fiber of $M\times_{B}\tilde{M} \into M$.
The space $M \times_B \tilde{M}$ is the {\it correspondence space}. 
\end{definition}

We have relaxed the definition one normally finds in the literature, where two integrality conditions are imposed: one requires that the closed 3-forms $H$ and  $\tilde{H}$ represent integral cohomology classes and that $F$ is {\it unimodular} in the sense that
\begin{equation}\label{eq:unimodular}
(F(\del_{\theta_i},\del_{\tilde{\theta_j}})) \in \Gl(k,\Z),
\end{equation}
where $\{\del_{\theta_i}\}$ is a basis of invariant period-1 elements for $\frak{t}_{M}^k$ and $\{\del_{\tilde{\theta_j}}\}$ is  a basis of invariant  period-1 elements for  $\frak{t}_{\tilde{M}}^k$.  These conditions are useful for interpreting the T-dual bundles as fibrations of mutually dual tori.
We now recall the salient implications of the T-duality relation, as pointed out in \cite{BEM03,BHM03,BuSc04}.

The tangent space of a $T^k$-bundle $M\into B$ fits in an exact sequence
$$ 0 \into \frak{t}^k \into TM/T^k \into TB,$$
and hence the space of invariant forms on $M$ is filtered:
$$ \Omega^k (B) = \mc{F}^0  \subset \mc{F}^1 \subset \cdots \subset \mc{F}^k = \Omega_{T^k}(M),$$ 
where $\mc{F}^i = \Ann(\wedge^{i+1} \frak{t})$.

The first implication of Definition \ref{def:t-dual} is that if $M$ and $\tilde{M}$ are T-dual, then $H \in \mc{F}^1$, i. e.,
\begin{equation}\label{eq:zero holonomy}
H(X,Y,\cdot)=0 \qquad \forall X,Y \in \frak{t}_M^k \subset TM.
\end{equation}
Indeed, this follows from \eqref{eq:T-duality condition} as $\tilde{p}^*\tilde{H}(X,\cdot,\cdot) = 0$ whenever $X \in \ker\tilde{p}_*$ and, since $F$ is invariant, $dF (X,Y,\cdot) = 0$ whenever $X$ and $Y$ are tangent to the $T^{2k}$-fiber.

The choice of a connection $\theta \in \Omega^1(M,\frak{t})$ gives a splitting for the exact sequence above, rendering $TM/T^k \cong \frak{t} \oplus TB$, and hence the filtration of invariant forms becomes a splitting. For example,
$$\Omega^3_{T^k}(M) = \sum_{i=0}^3 \Omega^i(M,\wedge^{3-i}\frak{t}^*).$$
Condition \eqref{eq:zero holonomy}, means that for the 3-form $H$ two of the four terms in this sum vanish and we have
\begin{equation}\label{eq:H looks like}
H = \IP{\tilde{c},\theta} +h
\end{equation}
with $\tilde{c} \in \Omega^2(B,\frak{t}^*)$ and $h \in \Omega^3(B).$ Since $H$ is closed, so is $\tilde{c}$. Further, if $H$ represents an integral cohomology class, so does $\tilde{c}$.

With such a choice of connection, \cite{BHM03} establishes that if $H$ represents an integral cohomology class, then \eqref{eq:zero holonomy} is the only obstruction to the existence of T-duals.

\begin{proposition}[Bouwknegt--Hannabuss--Mathai \cite{BHM03}]\label{prop:existence of T-duals}
Given a principal $k$-torus bundle $M$ over a base manifold $B$ and a 3-form $H \in \Omega^3_{T^k}(M)$ representing an integral cohomology class satisfying \eqref{eq:zero holonomy}, there is another $T^k$-bundle $\tilde{M}$ over $B$ T-dual to $M$.
\end{proposition}
\begin{proof}
Choose a connection $\theta \in \Omega^1(M,\frak{t})$ as above so that $H$ is given by \eqref{eq:H looks like} with  $\tilde{c}$ representing an integral cohomology class in $H^2(B,\frak{t}^{k*})$.
 
Now, we regard $\frak{t}^{k*}$ as the Lie algebra for the dual $k$-torus and let $\tilde{M}$ be a principal torus bundle over $B$ with Chern class $[\tilde{c}]$. Let  $\tilde{\theta} \in \Omega^1(\tilde{M},\frak{t}^{k*})$ be a connection on $\tilde{M}$ such that $d\tilde{\theta} = \tilde{c}$ and define
$$ \tilde{H} = \IP{c,\tilde{\theta}} + h \in \Omega_{cl}^3(\tilde{M}),$$
where $c = d\theta \in \Omega^2(M,\frak{t}^k)$. We claim that $(M,H)$ and $(\tilde{M},\tilde{H})$ are T-dual spaces. Indeed, according to our choices, on the correspondence space, we have
$$p^*H -\tilde{p}^*\tilde{H} = \IP{\tilde{c},\theta} - \IP{c,\tilde{\theta}}  = \IP{d\tilde{\theta},\theta} - \IP{d\theta,\tilde{\theta}}  = -d\IP{\theta,\tilde{\theta}}$$
So,  $(M,H)$ and $(\tilde{M},\tilde{H})$ are T-dual, and with these choices,  $F = - \IP{\theta,\tilde{\theta}}$ is unimodular.
\end{proof}

\begin{example}\label{ex:S1 duality}
A particular case of the construction given in Proposition \ref{prop:existence of T-duals} which will appear frequently concerns the case of circle bundles. In this case, the skew symmetry of $H$ implies that \eqref{eq:zero holonomy} holds trivially and hence any $S^1$-bundle has a T-dual. According to the construction above, if $(M,H)$ is a circle bundle with a closed 3-form representing an integral class, once we choose an identification $\frak{s}^1 \cong \R$, one T-dual can be constructed by choosing a connection $\theta \in \Omega^1(M)$ and using it to  write
$$H = \tilde{c} \wedge \theta + h,$$
with $\tilde{c}$ and $h$ basic forms. So a circle bundle $\tilde{M}$ with Chern class $[\tilde{c}]$, connection $\tilde{\theta}$ with curvature $d\tilde{\theta} =\tilde{c}$ and corresponding 3-form
$$\tilde{H} = c \wedge \tilde{\theta} + h$$
is T-dual to $(M,H)$.
\end{example}

Some concrete examples of the situation depicted above are as follows.
\begin{example}\label{ex:s3}
The Hopf fibration exhibits the 3-sphere, $S^3$, as a principal $S^1$ bundle over $S^2$. The curvature of this bundle is a volume form of $S^2$, $\sigma$. So $S^3$ with zero 3-form is $T$-dual to $(S^2 \times S^1, \sigma \wedge \tilde{\theta})$, where $\tilde{\theta}$ is an invariant volume for for $S^1$. 

On the other hand, if we consider the 3-sphere endowed with the 3-form $H=\theta \wedge \sigma$ (where $\theta$ is a connection 1-form such that $d\theta = \sigma$)  as a Hopf fibration over $S^2$, then it is T-dual to itself.
\end{example}

\begin{example}[Lie groups]\label{ex:lie groups} Let $(G,H)$ be a semi-simple Lie group with its Cartan 3-form $H(X,Y,Z) = K([X,Y],Z)$, where $K$ is the Killing form.  Recall that $[H]$ generates $H^3(G,\Z)$.

Choosing a subgroup $S^1 \subset G$ acting by left translations, we may view $G$ as a principal circle bundle. For $\del_{\theta}\in \frak{g}$ tangent to $S^1$ and of Killing length $-1$, a natural connection on $G$ is given by $-K(\del_{\theta}, \cdot)$. The curvature of this connection is given by
$$c(Y,Z) = d(-K(\del_{\theta}, \cdot))(Y,Z) = K(\del_{\theta}, [Y,Z])= H(\del_{\theta},Y,Z),$$
hence $c$ and $\tilde{c}$ are related by
$$c =  i_{\del_{\theta}} H = \tilde{c}.$$
This shows that semi-simple Lie groups equipped with their Cartan 3-forms are self T-dual.
\end{example}

T-duals are by no means unique, i.e., knowing the torus bundle $(M,H)$ does not determine $(\tilde{M},\tilde{H})$. This should be clear already from the argument used in Proposition \ref{prop:existence of T-duals}, as $\tilde{M}$ was chosen to be  `a torus bundle with Chern class $[\tilde{c}]$'. Since $[\tilde{c}]$ is simply the real image of an integral class $H^2(B,\Z) \into H^2(B,\R)$, the torsion part of the Chern class of $\tilde{M}$  is not determined by $\tilde{c}$ and can be chosen arbitrarily. Different choices for torsion will render different bundles. Furthermore, if we choose a basis making $\frak{t}^* \cong \R^k$ and $\tilde{c} = (\tilde{c}_1,\cdots,\tilde{c}_k)$ with $\tilde{c}_i \in \Omega_{cl}^2(B)$, one may modify the  construction from Proposition \ref{prop:existence of T-duals} so that a $T^k$-bundle with Chern class $(m_1[\tilde{c}_1],\cdots,m_k[\tilde{c}_k])$, $m_i \in \Z$, is also T-dual to $(M,H)$ (however in this case $F$ is no longer unimodular). 

The above sources of nonuniqueness are related to the fact that we are working over $\R$, while the Chern classes lie in $H^2(B,\Z)$. However nonuniqueness of T-duals goes beyond torsion.
The relevant remark is that the cohomology class of the component of $H$ lying in $\Omega^2(B,\frak{t}^*)$, called $[\tilde{c}]$,  is only defined in $H^2(B,\frak{t}^*)$ modulo $d( \Omega^0(B,\wedge^2\frak{t}^*))$.  Therefore, for $k>1$ the class $[\tilde{c}] \in H^2(B,\frak{t}^*)$ is not well defined. We illustrate this phenomenon in a concrete example. 

\begin{example}[Bunke--Schick \cite{BuSc04}] 
Let $M$ be a principal 2-torus bundle over $B$ with $H=0$. Choose  an identification $\frak{t} \cong \R^2$ and a connection $(\theta_1,\theta_2)$ for $M$ so that $c_i = d\theta_i$ represent the Chern classes of $M$. Following the construction from in Proposition \ref{prop:existence of T-duals}, we see that  the trivial 2-torus bundle $B \times  S^1 \times S^1$ with 3-form
$$\tilde{H} = c_1 \wedge \tilde{\theta_1}  + c_2 \wedge \tilde{\theta_2}.$$
is T-dual to $M$.

A less obvious T-dual is the bundle $\tilde{M}$ over $B$ with Chern classes $[c_1]$ and $-[c_2]$ (equivalent to changing $H =0$ by the exact 3-form $d(\theta_1\wedge \theta_2)$ and then using the construction of Proposition \ref{prop:existence of T-duals}). We can see that $\tilde{M}$ is T-dual to $M$ by choosing a connection $\tilde{\theta} \in \Omega^1(\tilde{M},\R^2)$ so that $d\tilde{\theta} =(c_1, -c_2)$ and choosing 
$$\tilde{H} = d(\tilde{\theta_1} \wedge \tilde{\theta_2}) = c_1 \wedge \tilde{\theta_2} +c_2\wedge\tilde{\theta_1}.$$ 
For this $\tilde{H}$ we have
\begin{align*}
p^*(H) - \tilde{p}^*(\tilde{H}) &= d(\theta_1)\wedge \tilde{\theta_2} + (d\theta_2) \wedge \tilde{\theta_1}= d(\theta_1\wedge \tilde{\theta_2} + \theta_2 \wedge \tilde{\theta_1}) - (\theta_1 \wedge c_2 - \theta_2 \wedge c_1)  \\& =  d(\theta_1\wedge \tilde{\theta_2} + \theta_2 \wedge \tilde{\theta_1}+\theta_1 \wedge \theta_2).
\end{align*}
So, $F = \theta_1\wedge \tilde{\theta_2} + \theta_2 \wedge \tilde{\theta_1}+\theta_1 \wedge \theta_2$ is unimodular and makes  $(M,\tilde{M})$ a  T-dual pair. 
\end{example}

Having described the basic behaviour of the T-duality relation between principal torus bundles with 3-form flux, we proceed to the consequence of this relation which principally concerns us. 
\begin{theorem}[Bouwknegt--Evslin--Mathai \cite{BEM03}]\label{theo:BEM}
If $(M,H)$ and $(\tilde{M},\tilde{H})$ are T-dual, with $p^*H - \tilde{p}^*\tilde{H} = dF$, then the following map
\begin{equation}\label{eq:tau}
\tau: \gO^{\bullet}_{T^k}(M) \lra \gO^{\bullet}_{T^k}(\tilde{M}) \qquad \tau(\rho) = \int_{T^k} e^{F}\wedge \rho,
\end{equation}
is an isomorphism of the differential complexes $(\Omega^{\bullet}_{T^k}(M),d_H)$ and $(\gO^{\bullet}_{T^k}(\tilde{M}),d_{\tilde{H}})$, where the integration above is along the fibers of $M \times_B \tilde{M} \into \tilde{M}$.
\end{theorem}
\begin{proof}[Sketch of the proof]
The proof consists of two parts. The first, which is just linear algebra, consists of proving that $\tau$ is invertible. This boils down to the requirement that $F$ satisfies \eqref{eq:nondegenerate}. The second part is to show that $\tau$ is compatible with the differentials, which is a consequence of equation \eqref{eq:T-duality condition}:
\begin{align*}
d_{\tilde{H}}\tau(\rho) &=\int_{T^k} d_{\tilde{H}}(e^{F} \wedge \rho)\\
&= \int_{T^k} (H-\tilde{H})\wedge e^{F} \wedge\rho + e^{F} d\wedge  \rho + \tilde{H}\wedge e^{F} \wedge \rho\\
&= \int_{T^k}  H \wedge e^{F} \wedge \rho + e^{F}\wedge  d \rho\\
& = \tau(d_H \rho).
\end{align*}
\end{proof}

\section{T-duality as a map of Courant algebroids}\label{courant algebroids}

In this section we rephrase the  T-duality relation as an isomorphism of Courant algebroids.  This will allow us to transport any $T^k$-invariant generalized geometrical structures on  $M$ to similar structures on $\tilde M$, and vice versa.

The map $\tau$ introduced in Theorem \ref{theo:BEM} may be described as  a composition of pull-back, $B$-field transform and push-forward, all operations on the Clifford module of $T^k$-invariant differential forms.  To make $\tau$ into an isomorphism of Clifford modules, we need to specify an isomorphism of $T^k$-invariant sections $\gf:(TM \oplus T^*M)/T^k \lra  (T\tilde{M} \oplus T^*\tilde{M})/T^k$ such that
\begin{equation}\label{eq:gf and tau}
\tau(v \cdot \rho) = \gf(v)\cdot \tau(\rho),
\end{equation}
for all $ v \in \Gamma(TM \oplus T^*M)/T^k$ and  $\rho \in \Omega^{\bullet}_{T^k}(M)$.   We now define such a map $\gf$.

Given T-dual spaces $(M,H)$ and $(\tilde{M},\tilde{H})$, consider the diagram
$$\xymatrix{ & (M\times_{B}\tilde{M}, p^*H - \tilde{p}^*
\tilde{H})\ar[dl]_p\ar[dr]^{\tilde{p}} &
\\ (M,H)\ar[dr]_\pi& &(\tilde{M},\tilde{H})\ar[dl]^{\tilde{\pi}} \\ &B & \\}
$$
Given $X + \xi \in (TM \oplus T^*M)/T^k$, we can try and pull it back to the correspondence space $M\times_B\tilde{M}$. While $p^*\xi$ is well defined, the same is not true about the vector part. If we pick a lift  of $X$, $\hat{X} \in T(M\times_B\tilde{M})$, any other lift differs from $\hat{X}$ by a vector  tangent to $T_{\tilde{M}}^k$, the fiber of the projection $p$. Then we form the $B$-field transform of $\hat{X} + p^*\xi$ by $F$ to obtain
\begin{equation}\label{eq:middle step}
\hat{X}  + p^*\xi -F(\hat{X}).
\end{equation}
We want to define $\gf$ as the push forward of the element above to an element of $(T\tilde{M}\oplus T^*\tilde{M})/T^k$, however there are two problems. First, $\tilde{p}_*(\hat{X})$ is not well defined, as it depends of the choice of lift $\hat{X}$ of $X$. Second, we can only push-forward the form $p^*\xi -F(\hat{X})$ if it is basic, i. e., only if
\begin{equation}\label{eq:push forward}
\xi(Y) -F(\hat{X},Y) =0, \qquad \forall Y \in \frak{t}^k_M.
\end{equation}
As is often said, two problems are better than one: the nondegeneracy of $F$ on $T_M^k \times T^k_{\tilde{M}}$ means that there is only one lift $\hat{X}$ for which \eqref{eq:push forward} holds and hence we define $\gf$ as the push forward of \eqref{eq:middle step} for that choice of lift of $X$:
$$\gf(X + \xi) = \tilde{p}_*(\hat{X}) +  p^*\xi -F(\hat{X}).$$
If one traces the steps above in the definition of $\gf$, it is clear that it satisfies \eqref{eq:gf and tau}.
Further, the compatibility of $\tau$ with the differentials translates to a compatibility of $\gf$ with the differential structure of $M$ and $\tilde{M}$, as we now explain.

\begin{theorem}\label{theo:gf} Let $(M,H)$ and $(\tilde{M},\tilde{H})$ be T-dual spaces. The map $\gf$ defined above is an isomorphism of Courant algebroids, i. e., for $v_1,v_2 \in (TM \oplus T^*M)/T^k$
$$\IP{v_1,v_2} = \IP{\gf(v_1),\gf(v_2)} \qquad \mbox{and} \qquad \gf([v_1,v_2]_H) = [\gf(v_1),\gf(v_2)]_{\tilde{H}}.$$
\end{theorem}

\begin{proof}
Both of these properties follow from \eqref{eq:gf and tau}. To show that $\gf$ is orthogonal, let $v \in TM \oplus T^*M$  and $\rho \in \Omega_{T^k}(M)$ and compute
$$\IP{v,v}\tau(\rho) = \tau(\IP{v,v} \rho) = \tau(v\cdot(v\cdot \rho)) = \gf(v)\cdot(\gf(v) \cdot \tau(\rho)) = \IP{\gf(v),\gf(v)} \tau(\rho).$$
To prove compatibility with the brackets we use equation \eqref{eq:twisted bracket} together with Theorem \ref{theo:BEM}:
\begin{align*}
\gf([v_1,v_2]_H)\cdot \tau(\rho) &  = \tau([v_1,v_2]_H \cdot\rho)\\
& = \tau([[d_H,v_1],v_2]\cdot\rho) \\
&= [[d_{\tilde H}, \gf(v_1)],\gf(v_2)]\cdot \tau(\rho)\\
&=[\gf(v_1),\gf(v_2)]_{\tilde H}\cdot \tau(\rho).
\end{align*}
\end{proof}

\begin{example}
In this example we give a concrete expression for the map $\gf$ introduced above in the case of the  T-dual $S^1$-bundles $(M,H)$ and $(\tilde{M},\tilde{H})$ constructed in Example \ref{ex:S1 duality}. Recall that $M$ and $\tilde{M}$ were  endowed with connections $\theta$ and $\tilde{\theta}$ for which $F = -\theta\wedge \tilde{\theta}$.

The presence of the connections $\theta$ and  $\tilde{\theta}$ provides us with a splitting $TM/S^1 \cong TB \oplus \IP{\del_\theta}$, where $TB$ is the space of invariant horizontal vector fields and $\del_{\theta}$ is an invariant  period-1 generator of the circle action. Similarly $T^*M/S^1 \cong T^*B \oplus \IP{\theta}$ and the same holds for $\tilde{M}$. So, an element in $(TM \oplus T^*M)/S^1$ can be written as
$$ X + f \del_{\theta} + \xi + g \theta,$$
with $X$ horizontal and $\xi$ basic. The pull back of this element to $M \times_B\tilde{M}$ is given by
$$ X + f \del_{\theta} + k \del_{\tilde{\theta}}+ \xi + g \theta,$$
where $k$ will be determined later. Then the $B$-field transform by $F$ is
$$ X + f \del_{\theta} + k \del_{\tilde{\theta}}+ \xi + g \theta +f\tilde{\theta} -k\theta.$$
The requirment that $ \xi + g \theta +f\tilde{\theta} -k\theta$ is basic for $M\times_B\tilde{M} \into \tilde{M}$ is equivalent to  $k =g$, and $\gf$ is defined as the push-forward of this element to $\tilde{M}$, yielding
\begin{equation}\label{eq:S1 duality map}
 \gf(X + f \del_{\theta} + \xi + g \theta) = X + g \del_{\tilde{\theta}}+ \xi  +f\tilde{\theta},
 \end{equation}
In the final expression, we recognize the exchange of tangent and cotangent  `directions' described by physicists in the context of  T-duality. 
\end{example}

\begin{remark}
Let $M$ and $\tilde{M}$ be T-dual spaces as  constructed in Proposition \ref{prop:existence of T-duals}. The choice of connections for $M$ and $\tilde{M}$ allows us to split the invariant tangent and cotangent bundles:
$$(TM\oplus T^*M)/{T^k} \cong TB\oplus T^*B \oplus \frak{t}^k \oplus \frak{t}^{k*},$$
$$(T\tilde{M}\oplus T^*\tilde{M})/{T^k} \cong TB\oplus T^*B \oplus \frak{t}^k \oplus \frak{t}^{k*}.$$
In this light, the map \gf\ is a the permutation of the terms $\frak{t}^k$ and $\frak{t}^{k*}$. This is Ben-Bassat's starting point for the study of T-duality and \gcss\ in \cite{BeB06,BeB06b}, where he deals with the case of linear torus bundles with vanishing background 3-form $H$.
\end{remark}

\section{T-duality and generalized structures}\label{transport}
If $(M,H)$ and $(\tilde M, \tilde H)$ form a T-dual pair, then we obtain a Courant isomorphism $\gf$ from Theorem \ref{theo:gf} between the $T^k$-invariant sections of $T\oplus T^*$ of $M$ and $\tilde M$.  This immediately implies that any $T^k$-invariant generalized geometrical structure, since it is defined purely in terms of the structure of $T\oplus T^*$, may be transported from one side of the T-duality to the other.  What is interesting is that the resulting structure may have very different behaviour from the original one, stemming from the fact that tangent and cotangent directions have been exchanged in the T-duality.  We explore this phenomenon in detail in this section.

\begin{theorem}\label{theo:gcss and t-duality}
Let $(M,H)$ and $(\tilde{M},\tilde{H})$ be T-dual spaces. Then any Dirac, generalized complex, \gk\ or SKT structure on $M$ which is invariant under the torus action is transformed via $\gf$ into a structure of the same kind on $\tilde{M}$.
\end{theorem}
Note also that the Courant isomorphism $\gf$ came together with a Clifford module isomorphism $\tau$  from Theorem \ref{theo:BEM}; it follows similarly that any decomposition of forms induced by generalized geometrical structures is preserved by $\tau$:  
\begin{cor}
Let $(M,H)$ and $(\tilde{M},\tilde{H})$ be T-dual spaces, and let $\J$ and $\tilde{\J}$ be a pair of T-dual \gcss\ on these spaces, as obtained from Theorem \ref{theo:gcss and t-duality}. Let $\mc{U}_M^k$ and $\mc{U}_{\tilde{M}}^k$ be the decompositions of forms induced by $\J$ and $\tilde{\J}$.   Then $\tau(\mc{U}^k_M) = \mc{U}^k_{\tilde{M}}$ and also
$$ \tau(\del_{M} \psi) = \del_{\tilde{M}}\tau(\psi) \qquad  \tau(\delbar_{M} \psi) = \delbar_{\tilde{M}}\tau(\psi),$$
for $\del_M,\del_{\tilde M}$ the generalized Dolbeault operators associated to $\J,\tilde\J$.
\end{cor}
A special case of the above occurs when the generalized complex structure $\J$ has a holomorphically trivial canonical bundle, in the sense that it has a non-vanishing $d_H$-closed section $\rho\in\Gamma(K)$. This additional structure is called a generalized Calabi-Yau structure~\cite{Hi03}.  In this case, any T-dual generalized complex structure $\tilde\J$ also has this property, with preferred generalized Calabi-Yau structure $\tilde\rho = \tau(\rho)$.   

In the following example we see how the behaviour of generalized complex structures may change under T-duality.  This change of geometrical type is at the heart of mirror symmetry, where complex and symplectic structures are exchanged between mirror Calabi-Yau manifolds.
\begin{example}[Change of type of \gcss]\label{ex:type change} Recall that the type of a \gcs\ $\rho = e^{B+i\omega} \wedge \Omega$ is the degree of the decomposable form $\Omega$.

Using this description, we determine how the types of  T-dual \gcss\ are related. If $\rho = e^{B+ i \omega}\wedge \Omega$ is a locally-defined invariant form defining a $T^k$ invariant \gcs\ on $M$,  then the corresponding \gcs\ on the T-dual $\tilde{M}$ is determined by the pure spinor $\tau(\rho)$, where $\tau$ is the map defined in \eqref{eq:tau}. Expanding the exponential, we see that the lowest degree form in $\tau(\rho)$ is determined by the smallest  $j$ for which
\begin{equation}\label{eq:type}
\int_{T^k}(F+B+ i\omega)^j \wedge\Omega \neq 0
\end{equation}
and hence the type of the T-dual structure $\tilde{\J}$ is
$$type(\tilde{\J}) = type(\J) +  2j -k,$$
where $j$ is the smallest  natural number for which \eqref{eq:type} holds.

For example, when performing T-duality along a circle, $j$ is either 1 or 0, depending on whether $\Omega$ is basic or not, respectively. Hence, the type of $\tilde{\J}$ is either $type(\J) +1$ or $type({\J})-1$, depending on whether  $\Omega$ is basic or not.

In order to obtain more concrete expressions for the relation between the type of T-dual structures, we make some assumptions about the algebraic form of one of the \gcss\ involved and of the form $F \in \Omega^2(M\times_B\tilde{M})$ which defines the T-duality.

If  $M^{2n}$ is a principal $T^n$ bundle with connection and $F = - \IP{\theta,\tilde{\theta}}$ is the form used in Proposition \ref{prop:existence of T-duals}, then we can determine the change in type for certain \gcss\ under T-duality, assuming the fibers are $n$-tori with special geometric constraints:
\begin{center}
\vskip6pt
\noindent
\begin{tabular}{|l|l|l|l|l|l| }
\hline

Structure on $M$ & Fibers of $M$ &$l$& $r$& Structure on $\tilde{M}$ & Fibers of $\tilde{M}$\\[3pt]
\hline
Complex & Complex & $n/2$ & $0$ & Complex & Complex \\
Complex & Real ($TF \cap I(TF) = \{0\}$)& $n$&  $0$ & Symplectic & Lagrangian\\
Symplectic & Symplectic & 0 & $n$ & Symplectic & Symplectic\\
Symplectic & Lagrangian & 0 & 0& Complex & Real\\
\hline
\end{tabular}
\vskip3pt
{\small Table 1: Change of type of \gcss\ under T-duality according to the type of fiber $F$.}
\vskip6pt
\end{center}
\end{example}

\begin{example}[Hopf surfaces]\label{ex:hopf surfaces}
Given two complex numbers $a_1$ and $a_2$, with $|a_1|, |a_2|>1$, the quotient of $\C^2$ by the action $(z_1,z_2) \mapsto (a_1 z_1,a_2 z_2)$ is a primary Hopf surface (with the induced complex structure). Of all primary Hopf surfaces, these are the only ones admiting a $T^2$ action preserving the complex structure (see \cite{BPVdV84}). If $a_1= a_2$, the orbits of the 2-torus action are elliptic surfaces and hence, according to Example \ref{ex:type change}, the T-dual will still be a complex manifold. If $a_1 \neq a_2$, then the orbits of the torus action are real except for the orbits passing through $(1,0)$ and $(0,1)$, which are elliptic. In this case, the T-dual will be generically symplectic except for the two special fibers corresponding to the elliptic curves, where there is type change. This example also shows that even if the initial structure on $M$ has constant type, the same does not need to be true in the T-dual.
\end{example}

\begin{example}[Mirror symmetry of Betti numbers]\label{ex:mirror symmetry and the Upq}
Consider the case of the mirror of a Calabi-Yau manifold along a
special Lagrangian fibration. We have seen that the bundles
$U_{\omega,I}^k$ induced by both the complex and symplectic structure
are preserved by T-duality. Hence $U^{p,q}=U_{\omega}^p\cap U_I^q$ is
also preserved, but, $U^{p,q}$ will be associated in the mirror to
$U_{\tilde{I}}^p\cap U_{\tilde{\omega}}^q$, as complex and symplectic
structure get swaped. Finally, as remarked in the previous section (equation \eqref{eq:upqkahler}), we have an isomorphism between
$\gO^{p,q}$ and $\mc{U}^{n-p-q,p-q}$. Making these identifications,
we have
$$\gO^{p,q}(M) \cong \mc{U}^{n-p-q,p-q}(M) \cong  \tilde{\mc{U}}^{n-p-q,p-q}(\tilde{M}) \cong \gO^{n-p,q}(\tilde{M}).$$
Which, in cohomology, gives the usual `mirror symmetry' of the Hodge diamond.
\end{example}

Since the map $\gf$ from Theorem \ref{theo:gf} is an orthogonal isomorphism, we may transport an invariant  generalized metric $\G$ on $(TM\oplus T^*M)/T^k$ to a generalized metric $\tilde{\G} = \gf \circ \G \circ \gf^{-1}$ on $(T\tilde{M}\oplus T^*\tilde{M})/T^k$. If $C_{\pm}$ are the $\pm1$-eigenspaces of $\G$ then $\tilde{C}_{\pm} = \gf(C_{\pm})$ are the $\pm 1$-eigenspaces of $\tilde{\G}$. This is a complete description of the T-dual generalized metric from the point of view of $TM\oplus T^*M$.

However, we saw that a generalized metric may also be described in terms of a Riemannian metric $g$ on $M$ and a real 2-form $b\in \Omega^2(M,\R)$, so that $C_+$ is the graph of $b + g$.  It is natural to ask how $g$ and $b$ change under T-duality.  We shall describe the dual metric $\tilde{g}$ and 2-form $\tilde{b}$ in the case of principal circle bundles, as in Example \ref{ex:S1 duality}.
\begin{example}[T-duality of generalized metrics]\label{ex:buscher}
Let $(M,H)$ and $(\tilde{M},\tilde{H})$ be T-dual principal circle bundles as in Example \ref{ex:S1 duality}, with connections $\theta$ and $\tilde{\theta}$. Let $\G$ be a  generalized metric on $M$ invariant \wrt\ the circle action, so that the induced metric $g$ on $TM$ and 2-form $b$ are both $S^1$-invariant and hence can be written as
\begin{align*}
g &= g_{0} \theta \odot \theta + g_{1} \odot \theta + g_{2}\\
b &= b_1 \wedge \theta + b_2,
\end{align*}
where the $g_i$ and $b_i$ are basic forms of degree $i$.
So, the elements of $C_+$, the $+1$-eigenspace of $\G$,  are of the form:
$$X + f \del_{\theta} + (i_X g_2 + f g_1 + i_X  b_2 - f b_1) + (g_1(X) + f g_0 + b_1(X))\theta.$$
Applying \gf, the generic element of $\gf(C_+) = \tilde{C_+}$ is given by:
$$X + (g_1(X) + f g_0 + b_1(X)) \del_{\tilde{\theta}} + (i_X g_2 + f g_1 + i_X b_2 - f b_1) +  f \tilde{\theta}.$$
This is the graph of $\tilde{b} + \tilde{g}$ for $\tilde b,\tilde g$ given by:
\begin{equation} \label{busher rules}
\begin{aligned}
\tilde{g} &= \frac{1}{g_0} \tilde{\theta} \odot \tilde{\theta} - \frac{b_1}{g_0}\odot \tilde{\theta} + g_2 + \frac{b_1\odot b_1 - g_1 \odot g_1}{g_0}\\
\tilde{b} & = - \frac{g_1}{g_0} \wedge \tilde{\theta} + b_2 + \frac{g_1 \wedge b_1}{g_0}
\end{aligned}
\end{equation}
These equations are certainly not new:  they were encountered by physicists in their computations of the dual Riemannian metric for T-dual sigma-models \cite{Bu87,Bu88} and carry the name {\it Buscher rules}.  The metric along the $S^1$ fiber, $g_0$, is sent to $g_0^{-1}$, considered the hallmark of a T-duality transformation.  
\end{example}

As we saw in Theorem~\ref{theo:bihermitiangk}, a generalized K\"ahler structure on $M$ induces a bi-Hermitian geometry $(g,I_\pm)$.  We wish to understand how this bi-Hermitian structure varies under T-duality.  Again, we consider the simple case of T-dual $S^1$-bundles.
\begin{example}[T-duality of bi-Hermitian structure]
The choice of a generalized metric $\G$ gives us two orthogonal spaces
$$C_{\pm} = \{X + b(X, \cdot) \pm  g(X,\cdot): X\in TM\},$$
and the projections $\pi_{\pm}:C_{\pm} \into TM$ are isomorphisms. Hence, any endomorphism $A\in End(TM)$ induces endomorphisms $A_{\pm}$ on $C_{\pm}$. Using the map \gf\ we can transport this structure to a T-dual:

\[
\xymatrix{A_+ \in \End(C_+)\ar[d]_{\pi} \ar[r]^\gf &\tilde{A}_+ \in \End(\tilde{C}_+)\ar[d]_{\pi} \\
A\in \End(TM) & \tilde{A}_{\pm}\in \End(T\tilde{M})\\
A_- \in \End(C_-) \ar[u]^{\pi}\ar[r]^{\gf}&\tilde{A}_- \in \End(\tilde{C}_-)\ar[u]^{\pi}
}
\]

As we are using the generalized metric to transport $A$ and the maps $\pi_{\pm}$ and \gf\ are orthogonal, the properties shared by $A$ and $A_{\pm}$ will be metric properties, e.g., self-adjointness, skew-adjointness and orthogonality. In the \gk\ case, it is clear that if we transport $I_{\pm}$ via $C_{\pm}$, we obtain the corresponding bi-Hermitian structure of the induced \gk\ structure on the T-dual:
$$\tilde{I}_{\pm} = \tilde{\pi}_{\pm} \gf\pi^{-1}_{\pm} I_{\pm} (\tilde{\pi}_{\pm} \gf\pi^{-1}_{\pm} )^{-1}.$$

In the case of a metric connection, $\theta = g({\del_\theta},\cdot)/g({\del_\theta},{\del_\theta})$, we can give a concrete description of $\tilde{I}_{\pm}$. We start describing the maps $\tilde{\pi}_{\pm} \gf\pi^{-1}_{\pm}$. If $V$ is orthogonal do $\del_\theta$, then $g_1(V)=0$ and
\begin{align*}
 \tilde{\pi}_{\pm} \gf\pi^{-1}_{\pm}(V) &= \tilde{\pi}_{\pm}\gf(V +b_1(V) \theta + b_2(V) \pm g_2(V,\cdot)) = \tilde{\pi}_{\pm}(V +b_1(V) \frac{\del}{\del \tilde{\theta}} + b_2(V) \pm g_2(V,\cdot)) \\
&=  V +b_1(V)\frac{\del}{\del \tilde{\theta}}.
\end{align*}
And for $\del_{\theta}$ we have
$$ \tilde{\pi}_{\pm} \gf\pi^{-1}_{\pm}(\del_{\theta}) =  \tilde{\pi}_{\pm} \gf(\del/\del \theta +b_1 \pm(\frac{1}{g_0}\theta  + g_1)) = \tilde{\pi}_{\pm}(\frac{1}{g_0}\del/\del \tilde{\theta} +\tilde{\theta})) = \pm\frac{1}{g_0}\del_{\tilde{\theta}}.$$

If  $V_{\pm}$ is the orthogonal complement to span$\{\del_{\theta}, I_{\pm}\del_{\theta}\}$, we may describe $\tilde{I}_{\pm}$ as follows:
\begin{equation}\label{eq:jpm}
\tilde{I}_{\pm}w = \begin{cases}I_{\pm}w,& \mbox{ if } w \in V_{\pm}\\
\pm\frac{1}{g_0} I_{\pm} \del_{\theta}& \mbox{ if } w = \del_{\tilde{\theta}}\\
\mp g_0 \del_{\tilde{\theta}} & \mbox{ if } w = I_{\pm} \del_\theta\\
\end{cases}
\end{equation}
Therefore, if we identify $\del_{\theta}$ with $\del_{\tilde{\theta}}$ and their orthogonal complements with each other via $TB$, $\tilde{I}_+$ is essentially the same as $I_+$, but stretched in the directions of $\del_{\theta}$ and $I_+ \del_{\theta}$ by $g_0$, while $\tilde{I}_-$ is $I_-$ conjugated and stretched in those directions. In particular, $I_+$ and $\tilde{I}_+$ induce the same orientation while $\tilde{I}_-$ and $I_-$ induce opposite orientations.
\end{example}
\begin{example}[T-duality and the generalized K\"ahler structure of Lie groups]
\label{ex:lie groups2}Any compact semi-simple Lie group $G$, together with its Cartan $3$-form $H$, admits a \gk\ structure (\cite{Gu03}, Example 6.39). These structures are obtained using the bihermitian point of view: any pair of left and right invariant complex structures on the Lie group $I_l$ and $I_r$, orthogonal \wrt\ the Killing form, satisfy the hypotheses of Theorem \ref{theo:bihermitiangk}.  Any \gk\ structure obtained this way will be neither left nor right invariant since at any point it depends on both $I_l$ and $I_r$. However one can also show that $I_l$ and $I_r$ can be chosen to be bi-invariant under the action of a maximal torus, and hence the corresponding \gk\ structure will also have this invariance. In this case, according to Theorem \ref{theo:gcss and t-duality} and Example \ref{ex:lie groups}, T-duality furnishes other \gk\ structures on the Lie group.

If we chose $I_+= I_r$ and $I_-=I_l$, the computation above shows that T-duality will furnish a new structure on the Lie group coming from $I_r$ and $\tilde{I}_l$, where $\tilde{I}_l$ is still left invariant  but induces the opposite orientation of $I_l$. Of course we may also swap the roles of $I_{\pm}$, changing the right invariant complex structure while the left invariant structure is fixed.
\end{example}

\begin{example}
An SKT structure is normally defined as a Hermitian structure $(g,I)$ on a manifold $M$ for which $dd^c \omega= 0$. According to Theorem \ref{theo:gcss and t-duality}, if we endow $M$ with the closed 3-form $d^c\omega$, then any T-dual  to $(M,H)$ obtains an SKT structure. The T-dual SKT structure $(\tilde{g},\tilde{I})$ is given by combining the Buscher rules for $\tilde g$ and the transformation law given in the previous example for the complex structure $I_-$.
\end{example}

\section{Explicit Examples}

In this section we study some instructive explicit examples of T-duality.

\begin{example}[The symplectic 2-sphere]\label{ex:s2}
Consider the standard circle action on the 2-sphere $S^2$ fixing the north and south poles $N,S$. we may view $S^2\backslash \{N,S\}$ as a trivial circle bundle over the interval $(-1,1)$. Using coordinates $(t,\theta) \in (-1,1) \times (0,2\pi)$, the round metric is given by
$$ds^2 = (1-t^2) d\theta^2 + \frac{1}{1-t^2}dt^2.$$
Using the Busher rules with $b=0$, the T-dual metric (introducing a coordinate $\tilde\theta$ along the T-dual fiber) is
$$d\tilde{s}^2 = \frac{1}{1-t^2} d\tilde{\theta}^2 + \frac{1}{1-t^2}dt^2.$$
Observe that the fixed points give rise to circles of infinite radius at a finite distance. This metric is not complete.

Given an invariant symplectic structure on the sphere $\omega =
w(t) dt\wedge d\theta$, and any invariant $B$-field
$B= b(t) dt\wedge d\theta$, we T-dualize the generalized complex structure defined by the differential form  $e^{B+i\omega}=1 + (b(t) + i w(t)) dt \wedge d \theta$. The dual
structure is given by
$$\tau (e^{B+i\omega}) =  d\tilde{\theta} +(b(t) + i w(t))dt.$$
Note that this defines a complex structure, with complex coordinate 
$$z(t,\tilde\theta) =e^{i\tilde\theta-\int_{-1}^t (\omega(t')-ib(t'))dt'} $$ 
Therefore the T-dual is biholomorphic to an annulus with interior radius 1 and exterior radius $\exp(-\int_{S^2} \omega)$.


While this is an accurate picture of what T-duality does to a (twice-punctured) symplectic sphere, from the physical point of view this is incomplete. In order for the symplectic sphere and the complex annulus to describe the same  physics, one must deform the complex structure to $\C^*$ and endow the resulting space with a complex-valued function called the {\it superpotential}, making it a Landau--Ginzburg model \cite{Ho02}.
\end{example}

\begin{example}\label{ex:gibbons--hawking}{\it (Odd 4-dimensional structures and the Gibbons--Hawking Ansatz)}
A \gcy\ metric structure on $(M,H)$ is a \gks\  $(\J_1,\J_2)$ such that the canonical bundles of $\J_1$ and $\J_2$ both admit nowhere-vanishing $d_H$-closed sections $\rho_1,\rho_2$, with the volume normalization
\begin{equation}\label{eq:cycondition}
(\rho_1,\overline{\rho_1}) = (\rho_2,\overline{\rho_2}).
\end{equation}
This last condition is the generalization of the usual Monge-Amp\`ere equation for Calabi-Yau metrics.

The description of \gcy\ geometry in real dimension 4 can be
divided in two cases, according to whether the bi-Hermitian induced complex
structures $J_{\pm}$ determine the same orientation or not. If they
determine different orientations, the corresponding differential forms $\rho_1$ and $\rho_2$ are of odd degree  and $J_{\pm}$ commute (see \cite{Gu03}, remark
6.14). The real distributions $S_{\pm} = \{v
\in TM: J_+ v = \pm J_-v\}$ are integrable \cite{ApGu06}, yielding a pair of
transverse foliations for $M$. If we choose holomorphic coordinates
$(z_1,z_2)$ for $J_+$ respecting this decomposition, then $(z_1,
\overline{z_2})$ furnish holomorphic coordinates for $J_-$.

As the metric $g$ is of type (1,1) \wrt\ both $J_{\pm}$, it is of
the form
$$ g= g_{1\overline{1}}dz_1 d\overline{z_{1}}+g_{2\overline{2}}dz_2 d\overline{z_{2}}.$$
The graphs of the $i$-eigenspaces of $J_{\pm}$  via $b\pm
g$ coincide with the intersections $L_1\cap L_2$ and $L_1 \cap
\overline{L_2}$ of the $+i$-eigenspaces $L_1,L_2$ of $\J_1,\J_2$ respectively.  Hence we can recover $L_1$ and $L_2$ from $J_{\pm}, g$
and $b$.  In this case, the differential forms annihilating
$L_1, L_2$ and hence generating the canonical bundles for $\J_1,\J_2$ are
$$\rho_1 = e^{b + g_{2\overline{2}}dz_2\wedge d \overline{z_2}} \wedge f_1 dz_1,\ \ \ \ \ \rho_2 = e^{b + g_{1\overline{1}}dz_1\wedge d\overline{z_1}} \wedge f_2 dz_2.$$
The \gcy\ condition $d\rho_{1} = d\rho_2=0$ implies that $f_1= f_1(z_1)$ is a
holomorphic function on $z_1$ and $f_2 = f_2(z_2)$ a holomorphic function on
$z_2$, hence with a holomorphic change of coordinates, we have
$$\rho_1 = e^{b + g_{2\overline{2}}dz_2\wedge d\overline{z_2}} \wedge dz_1,$$
$$\rho_2 = e^{b + g_{1\overline{1}}dz_1\wedge d\overline{z_1}} \wedge dz_2.$$
After rescaling $\rho_2$, if necessary, the compatibility condition
\eqref{eq:cycondition} becomes
$g_{1\overline{1}} = g_{2\overline{2}}$, showing that the
metric is conformally flat. Call this conformal factor $V$. 

Finally, the integrability conditions, $d\rho_1 = d\rho_2 =0$, impose
$$db\wedge dz_i = dV \wedge *dz_i = (* dV) \wedge dz_i, \qquad i = 1,2,$$
where $*$ is the Euclidean Hodge star. Therefore,
\begin{equation}\label{eq:Vandb}
db = * dV,
\end{equation}
showing that the conformal factor $V$ is harmonic with respect to
the flat metric.

Suppose the structure described above is realized on the 4-manifold $S^1
\times \R^3 $ in a way which is invariant by the obvious $S^1$ action, and remove a collection of points in $\R^3 $ so as to allow poles of $V$. The invariance of $V$
implies it is a harmonic function on $\R^3 $.   choosing the flat connection 1-form $\theta$, we may write $b =b_1\wedge \theta + b_2$.  Then, equation \eqref{eq:Vandb} implies that $db_1= *_3 dV$ and $db_2 =0$.  According to the Busher rules, the T-dual metric is given by
$$ \tilde{g}= V(dx_1^2 + dx_2^2+dx_3^2) + \frac{1}{V}(\tilde{\theta} - b_1)^2;$$
$$\tilde{b} = b_2,$$
with $db_1 = *_3 dV$ and $b_2$ closed.  This is nothing  but a $B$-field
transform of the hyper-K\"ahler metric given by the Gibbons--Hawking
ansatz.  
\end{example}

\section{Reinterpretations of T-duality}

In this section we show that it is possible to rephrase the definition of T-duality in two different ways:  as a double quotient and  as a submanifold. While the first point of view makes clear that the Courant algebroids $(TM\oplus T^*M)/T^k$ and $(T\tilde{M}\oplus T^*\tilde{M})/T^k$ are isomorphic, the second likens T-duality to a Fourier--Mukai transform, much in the spirit of \cite{DOPW01}.
\subsection{T-duality as a quotient}

In this section we review the process of reduction of Courant algebroids from \cite{BCG05}, and interpret T-duality in this light. We refer to \cite{BCG05} for more details on reduction.

Given a manifold with closed 3-form $(\M,\Hh)$, any section $v \in \Gamma(T \oplus T^*)$ defines a natural infinitesimal symmetry of  the orthogonal and Courant structures on $T \oplus T^*$, via the Courant adjoint action $\mathrm{ad}_v(w) = [v,w]_\Hh$.  

The fact that this is an infinitesimal symmetry is a consequence of the following properties of the Courant bracket:
\begin{align*}
\mc{L}_{\pi(v)}\IP{w_1,w_2} &= \IP{[v, w_1]_\Hh,w_2} + \IP{w_1,[v, w_2]_\Hh}\\
[v,[w_1,w_2] _{\Hh}]_\Hh& = [[v, w_1]_\Hh,w_2]_\Hh + [w_1,[v , w_2]_\Hh]_\Hh.
\end{align*}
Using the decomposition $T \oplus T^*$, we may write $\mathrm{ad}_{X+\xi}$ in block matrix form
\begin{equation}\label{eq:action matrix form}
\mathrm{ad}_{X+\xi}(Y+\eta) = \begin{pmatrix}\mc{L}_X & 0 \\ d\xi - i_X {{\Hh}} & \mc{L}_X\end{pmatrix} \begin{pmatrix}Y \\ \eta\end{pmatrix}.
\end{equation}
It is clear from this expression that $\mathrm{ad}_{X+\xi}$ consists of an infinitesimal symmetry of $\M$ (the Lie derivative $\mc{L}_X$) together with a $B$-field transform.

\begin{definition}
Let $G$ be a connected Lie group acting on a manifold $\M$ and $\psi:\frak{g} \into \Gamma(T\M)$ be the corresponding Lie algebra map. A {\it lift} of the action of $G$ to  $T\M \oplus T^*\M$ is a bracket-preserving map $\Psi:\frak{g} \into T\M \oplus T^*\M$, such that the following diagram commutes
\[
\xymatrix{
\frak{g} \ar[r]^{\Id} \ar[d]_{\Psi}& \frak{g}\ar[d] _{\psi}\\
\Gamma(T\M \oplus T^*\M)
\ar[r]^{~~~~\pi_T}& \Gamma(T\M)}
\]
and the infinitesimal $\frak{g}$ action induced by $\Psi(\Gg)$ integrates to a $G$ action on $T\M \oplus T^*\M$.
\end{definition}
Under these  circumstances, $T\M \oplus T^*\M$ is an equivariant $G$-bundle.  According to \eqref{eq:action matrix form}, the action of $\Psi(\gamma) = X_{\gamma}+\xi_\gamma$ preseves the splitting of $T\M \oplus T^*\M$ \iff\
\begin{equation}\label{eq:H condition}
i_{X_\gamma} \Hh = d\xi_{\gamma},
\end{equation}
in which case, according to \eqref{eq:action matrix form}, the  $\frak{g}$ action integrates to the standard $G$ action on $T \M$ and $T^*\M$ obtained by differentiation.  In particular, if $\Psi$ preserves the  splitting $T\M \oplus T^*\M$, it must be a lifted action, and hence  \eqref{eq:H condition} provides an effective way to check whether a bracket preserving map $\Psi$ is a lifted action. Conversely, given a lift of an action of a compact Lie group, one can always average $T\M$ for that lifted action to obtain a new splitting of $T\M \oplus T^*\M$ which is preserved by the action \cite{BCG05}. So, for compact groups, the requirement that the splitting is preserved is not restrictive.

\begin{remark}As an aside, given a lifted action $\Psi:\frak{g} \into \Gamma(T{\M}\oplus T^*{\M})$  for which \eqref{eq:H condition} holds, define  $\Psi = X +\xi$ with $X \in \Gamma(T\mc{\M}; \frak{g}^*)$ and $\xi \in \Gamma(T^*\mc{\M}; \frak{g}^*)$. Then combined form $\Hh + \xi$ is an equivariant 3-form in the Cartan complex, and
\begin{equation}\label{eq:dgH}
d_\frak{g}(\Hh + \xi) = \IP{\Psi(\cdot), \Psi(\cdot)} \in \mathrm{Sym}^2\frak{g}^*,
\end{equation}
showing that the pairing $\IP{\Psi(\frak{g}), \Psi(\frak{g})}$  is constant over $\M$.  Hence, a lift of an action gives rise to a symmetric 2-form on $\frak{g}$.\end{remark}

If the $G$-action on $\M$ is free and proper and $\Psi:\Gg \into \Gamma(T\M \oplus T^*\M)$ is a lift of $\psi$, then the distribution $K \subset T\M \oplus T^*\M$ generated by  $\Psi(\Gg)$ is actually a smooth subbundle. Furthermore, $K$ is $G$-invariant, and its orthogonal complement (with respect to the natural pairing), $\Kperp$, is also $G$-invariant.  Then $\Gamma^G(\Kperp)$ is closed under the bracket and 
$\Gamma^G(K \cap \Kperp)$ is an ideal in $\Gamma^G(\Kperp)$. Therefore
$$\E_{red} = \left. \frac{\Kperp}{K \cap \Kperp}\right/G, $$
as a bundle over  $\M/G$, inherits a bracket as well as a nondegenerate pairing. This is the main argument in the following theorem which is a particular case of Theorem 3.3 from \cite{BCG05}.
\begin{theorem}
Let $G$ act freely and propertly on $\M$, and let  $\Psi:\frak{g} \into T\M \oplus T^*\M$ be a lift of the $G$-action satisfying \eqref{eq:H condition}. Letting  $K = \Psi(\frak{g})$, the distribution
$$\E_{red} = \left. \frac{\Kperp}{K \cap \Kperp}\right/G $$
is a bundle over $\M/G$ which inherits a bracket and a nondegenerate pairing, making it into a Courant algebroid over $\M/G$. This Courant algebroid is exact\footnote{An exact Courant algebroid over $N$ is one which is locally isomorphic to $TN\oplus T^*N$.} \iff\ $K$ is isotropic.
\end{theorem}

\begin{definition}
The bundle $\E_{red}$ defined over $\M_{red}= \M/G$, equipped with its induced bracket and pairing, is the {\it reduced Courant algebroid} and $\M_{red}$ is the {\it reduced manifold}.
\end{definition}

Observe that the condition that $K$ is isotropic is equivalent to the requirement that the symmetric pairing~\eqref{eq:dgH} vanishes.  This would imply that $\Hh+\xi$ defines an equivariantly closed extension of $\Hh$. 

We now specialize to two examples of this quotient construction which will be relevant for T-duality. We distinguish them according to the symmetric form induced on the Lie algebra.

\begin{example}[Isotropic actions]\label{ex:H equiv closed} As we have just seen, to obtain {\it exact} reduced Courant algebroids, one must require that the lifted action, $\Psi = X + \xi$, is isotropic, i.e. $\Hh + \xi$ must be equivariantly closed. In this example we study this setting in detail and provide an explicit isomorphism between the reduced algebroid and $T{\M}/G\oplus T^*{\M}/G$.

In order to fully describe the reduced algebroid over $\M/G$, we choose a connection $\theta \in \Omega^1(\M;\Gg)$ for $\M$, viewed as a principal $G$ bundle. The extended action is given by $\Psi = X + \xi$ with $X \in \Gamma(T\M; \frak{g}^*)$ and $\xi \in \Gamma(T^*\M; \frak{g}^*)$.  We then apply a B-field transform by $B = \IP{\theta,\xi}  +\frac{1}{2}\IP{\theta\wedge\theta,\IP{\xi,X}}\in \Omega^{2}(\M,\R)$ to $T\M \oplus T^*\M$, so that the generators of the action become
$$X + \xi - i_X\IP{\theta,\xi} - i_X  \frac{1}{2}\IP{\theta\wedge\theta,\IP{\xi,X}}= X+\xi - \xi  +\IP{\theta,\IP{\xi,X}} -2 \frac{1}{2}\IP{\theta,\IP{\xi,X}} = X$$
and the 3-form curvature of $T\M \oplus T^*\M$ becomes the basic 3-form $\Hh_{red} = {\Hh} + dB$. After this $B$-field transform, the lifted action lies in $T\M$, and hence $\Kperp  = T\M + \Ann(\psi(\frak{g}))$, and
$$\left.\frac{\Kperp}{K}\right/ G \cong \left.(T\M/\psi(\frak{g}) \oplus \Ann(\psi(\frak{g}))\right/ G \cong T\M/G \oplus T^*\M/G,$$
where $\M/G$ is endowed with the $3$-form $\Hh_{red}$ defined above.

\end{example}

\begin{example}\label{ex:symplectic K} In this example, consider a lifted action of a Lie group $\mc{G}$ of dimension $2k$, with Lie algebra $\frak{G}$, for which the pairing induced on $\frak{G}$ is nondegenerate with split signature.  As we will show shortly, this is the case which arises in T-duality.

As before, let $\mc{K} = \Psi(\frak{G})$.  The nondegeneracy of the pairing on $\mc{K}$ implies that $\mc{K}^\perp \cap \mc{K} = \{0\}$, hence the reduced Courant algebroid is not exact.  Also, it is given by the $\mc{G}$-invariant sections of $\mc{K}^\perp$:
$$\E_{red} = \left.\frac{\mc{K}^\perp}{\mc{K} \cap \mc{K}^\perp}\right/\mc{G} = \mc{K}^\perp/\mc{G}.$$

Since the natural pairing has split signature on $\mc{K}$, we can choose a polarization of $\mc{K}$, expressing it as a sum of two isotropic subspaces, $\mc{K} = K\oplus \tilde{K}$. We show that such a decomposition provides us with alternative descriptions of  $\E_{red}$. Namely, it is clear that $K^\perp = \mc{K}^\perp + K$ and similarly for $\tilde{K}$ hence, {\it as vector bundles with a symmetric pairing},
\begin{equation}\label{eq:reduced E}
\E_{red} \cong \left.\frac{K^\perp}{K}\right/\mc{G} \cong \left.\frac{\tilde{K}^{\perp}}{\tilde{K}}\right/\mc{G}.
\end{equation}
Note, however, that since we did not require $K$ and $\tilde{K}$ to be Courant integrable, the resulting descriptions of $\E_{red}$ do not inherit natural Courant brackets.

We remedy this in the case where $\mc{G}$ is a product, $\mc{G} = G^k \times \tilde{G}^k$, and each of $\frak{g},\tilde{\frak{g}}$ is isotropic in $\frak{G} = \frak{g}\oplus \tilde{\frak{g}}$ for the induced pairing.  In this case, we are presented with a natural choice of isotropic splitting $\mc{K} = \Psi(\frak{g})\oplus \Psi(\tilde{\frak{g}}) =  K \oplus \tilde{K}$, and since both $K$ and $\tilde{K}$ are closed \wrt\ the Courant bracket, the spaces $(\Kperp/K)/(G\times \tilde{G})$ and $(\tilde{K}^{\perp}/\tilde{K})/(G \times \tilde{G})$ inherit a Courant bracket which agrees with that on $\E_{red} = \mc{K}^{\perp}$.  In this way, we obtain two alternative descriptions of $\E_{red}$.

In light of Example \ref{ex:H equiv closed}, the description $\E_{red} \cong (\tilde{K}^\perp/\tilde{K})/\tilde{G}\times G$ with $\tilde{K}$ isotropic  is very suggestive. Indeed, the space $(\tilde{K}^\perp/\tilde{K})/\tilde{G}$ is precisely the reduction of $T\mc{M}\oplus T^*\M$ by the isotropic action of $\tilde{G}$ to $M = \mc{M}/\tilde{G}$, obtaining an exact Courant algebroid over that space. Hence the reduced algebroid $\E_{red} \cong (\tilde{K}^\perp/\tilde{K})/\tilde{G}\times G$ corresponds to the quotient vector bundle $((\tilde{K}^\perp/\tilde{K})/\tilde{G})/ G$, which is isomorphic to the $G$-invariant sections of $TM \oplus T^*M$, viewed as a bundle  over $B= M/G$.

Reversing the roles of $G$ and $\tilde{G}$, we see that {\it the same Courant algebroid} can be described as the quotient vector bundle $((\tilde{K}^\perp/\tilde{K})/G)/ \tilde{G}$ which corresponds to the quotient bundle of the exact Courant algebroid $(\tilde{K}^\perp/\tilde{K})/G \cong T\tilde{M} \oplus T^*\tilde{M}$ over $\tilde{M} = \M/G$, so we have the following diagram:



\begin{equation}\label{eq:T-duality diagram3}
\xymatrix{ &T\M \oplus T^*\M           \ar[dl]_{\tilde{K}}           \ar[dr]^{K} &\\
(\tilde{K}^{\perp}/K)/\tilde{G}\ar[dr]_{/G} &            &(\Kperp/K)/G\ar[dl]^{/\tilde{G}}\\
&\mc{K}^{\perp}/G \times \tilde{G}&  \\}
\end{equation}
\end{example}

\begin{theorem}\label{theo:T-duality and reduction}
Let $(\M,\mc{H})$ be the total space of a principal $T^k \times  \tilde{T}^k$-torus bundle and let $\Psi:\frak{t}^k\times \frak{\tilde{t}}^k \into \Gamma(T\M\oplus T^*\M)$ be a lift of the $2k$-torus action for which the natural pairing on $\mc{K} = \Psi(\frak{t}^k\times \frak{\tilde{t}}^k)$  is nondegenerate, of split signature, and such that $\Psi(\frak{t}^k)$ and $\Psi(\tilde{\frak{t}}^k)$ are isotropic. Then the spaces $M,\tilde M$ obtained by reducing $\M$ by the action of $T^k$ and $\tilde{T}^k$ are T-dual.
Conversely, any pair of T-dual spaces arises in this fashion.
\end{theorem}
\begin{proof} By the argument from Example \ref{ex:H equiv closed}, one can transform $T\M \oplus T^*\M$ by a $B$-field so that  the lift of the action of $\tilde{T}^k$ lies in $T\M$. From now on, we assume that the action is of this form and denote by $\mc{H}$ the 3-form curvature of $T\M\oplus T^*\M$ associated to this splitting. Therefore, it follows from Example \ref{ex:H equiv closed} that the reduced space for the $\tilde{T}^k$ is $M = \M/\tilde{T}^k$ and that $\mc{H}$ is the pull back of $H$, the 3-form on $M$.

By the same argument,  one can transform $T\M \oplus T^*\M$ by a $B$-field $F$ so that  the lift of the action of $T^k$ lies in $T\M$. The reduction of $\M$ by this action gives $\tilde{M} = \M/T^k$ as a reduced space, with 3-form given by $\tilde{H} = \Hh + dF$.  So, $\M$ is the correspondence space for $M$ and $\tilde{M}$, and there we have $H -\tilde{H} = dF$.

Finally, since the pairing on $K$ is nondegenerate, the 2-form $F$ defined above gives rise to a nondegenerate pairing between $\psi(\frak{t}^k)$ and $\psi(\tilde{\frak{t}}^k)$, proving that $(M,H)$ and $(\tilde{M},\tilde{H})$ are T-dual.

This procedure can be reversed to prove the converse. Namely, given T-dual spaces $(M,H)$ and $(\tilde{M},\tilde{H})$,  let $\M$ be the correspondence space and define a lifted $T^k \times \tilde{T}^k$-action infinitesimally by
$$\Psi(t,\tilde{t}) = X_t - i_{X_t}F + X_{\tilde{t}},$$
where $X_t$ and $X_{\tilde{t}}$ are the vector fields associated to the Lie algebra elements $t$ and $\tilde{t}$.  The individual lifts of the actions of $T^k$ and $\tilde{T}^k$ are isotropic, but nondegeneracy of $F$ means that the natural pairing restricted to $\mc{K} = \Psi(T^k \times \tilde{T}^k)$ is nondegenerate, with split signature.
\end{proof}

\begin{remark}
From this point of view, Theorem \ref{theo:gf} may be viewed as the observation that the Courant algebroid of $T^k$-invariant sections of $TM\oplus T^* M$ is isomorphic to the algebroid of $\tilde{T}$-invariant sections of $T\tilde{M}\oplus T^*\tilde{M}$ because they are both isomorphic to the reduction of the correspondence space $\M = M\times_B\tilde M$  by the full $T^k \times \tilde{T}^k$ action.
\end{remark}

The approach to T-duality via quotients has also been described by Hu in \cite{Hu05b}, where he drives the theory much further, using it to study nonabelian Poisson group duality, essentially using this theorem as a definition of the duality relation.  There is also a strong similarity between this approach and the gauged sigma model approach in~\cite{RoVe92}.

\subsection{T-duality as a generalized submanifold}

In this section we re-interpret the T-duality relation between generalized complex manifolds $(M,H,\J)$ and $(\tilde M, \tilde H,\tilde\J)$ as a special kind of submanifold in the product $M\times\tilde M$.  

Given a manifold with 3-form $(\mc{N},\mc{H})$, a generalized  submanifold is a submanifold $\iota:\mc{M} \into \mc{N}$ together with a 2-form $F \in \Omega^2(\mc{M})$ such that $dF = \iota^*\mc{H}$. For any such generalized submanifold, one can form the generalized tangent bundle
$$\tau_F =\{X + \xi \in T\mc{M} \oplus T^*\mc{N}:\xi|_{\mc{M}} = F(X)\}.$$

If $(\mc{N},\mc{H})$ is endowed with a \gcs\ $\J$, we say that $(\mc{M},F)$ is a generalized complex submanifold if $\tau_F$ is invariant under $\J$.  When $\mc{N}$ is a usual complex manifold, this condition specializes to the usual notion of complex submanifold, and when $\mc{N}$ is symplectic, Lagrangian submanifolds provide examples.  

\begin{theorem}
Let $(M,H)$ and $(\tilde{M},\tilde{H})$ be two principal $k$-torus bundles over a manifold $B$. Let $(\mc{N},\mc{H}) = (M \times \tilde{M},H - \tilde{H})$ be their product and $\mc{M} = M \times_B \tilde{M}$  their fiber product, so that we have an inclusion $\iota:\M \into \mc{N}$. Then $M$ and $\tilde{M}$ are T-dual \iff\ there is a 2-form $F \in \Omega^2(\M)$ making $(\M,F)$ into a generalized submanifold and such that 
$\tau_F$ is everywhere transversal to
$$\tau_M = TM\oplus T^*M \subset T\mc{N} \oplus T^*\mc{N},$$
and similarly for $\tilde{M}$.
\end{theorem}

\begin{proof}
The requirement that $\mc{M}$ is a generalized submanifold is nothing but the condition
$$p^* H - \tilde{p}^*\tilde{H} = dF  $$
from the definition of T-duality.
Now consider
$$\tau_F \cap \tau_M =  \{X + \xi  \in TM \cap T\mc{M} \oplus T^*M : \xi|_{\M} = i_X F\}.$$
Let $X +\xi$ be an element of the above intersection. Since the projection $\mc{M} \into \tilde{M}$  has kernel the fibers $T_M^k$ and $\tilde{p}_* X = 0$,  we see that $X$ must be tangent to $T_M^k$. Therefore, requiring that $\tau_F \cap \tau_M = \{0\}$ is equivalent to requiring that $F(X) \in \Omega^1(\M)$ is not a pull back from $M$ for each $X \in T_M^k$, i.e. , that $F:T_M^k \times T_{\tilde{M}}^k \into \R$ is nondegenerate.
\end{proof}

As we now show, this viewpoint allows us to express the T-duality relation between generalized complex manifolds $(M,\J),(\tilde M, \tilde\J)$ as the existence of a generalized complex submanifold in the product $M\times\tilde M$.  In this way, the T-duality relation for generalized complex structures may be viewed as a generalization of the Fourier-Mukai transform for complex manifolds.

\begin{theorem}
Let   $(M,H)$ and $(\tilde{M},\tilde{H})$ be a T-dual pair and let $\J,\tilde{\J}$ be \gcss\ on these spaces. Endow the product $(\mc{N},\mc{H}) = (M\times \tilde{M},H - \tilde{H})$ with the \gcs\ $(\J , c \tilde{\J} c^{-1})$,
where $c:T\tilde{M}\oplus T^*\tilde{M} \into T\tilde{M} \oplus T^*\tilde{M}$ is given in matrix form by
$$\begin{pmatrix}
1& 0\\0 & -1
\end{pmatrix}.$$
Then $\J$ and $\tilde{\J}$ are T-dual \iff\ the correspondence space $(\M,F)$ is a generalized complex submanifold of $(N,\mc{H})$. 
\end{theorem}
\begin{proof}
Observe that the generalized tangent space of $(\M,F)$ is given by
$$\tau_F = \{(v,  c\gf (v)) \in (TM \oplus T^*M) \oplus (T\tilde{M}\oplus T^*\tilde{M})\}.$$ 
Hence this space is invariant under $(\J, c \tilde\J c^{-1})$ \iff\
$$c \J(v) = c\gf \tilde{\J}\gf^{-1}(v),\ \ \ \forall v\in TM\oplus T^*M,$$
that is, \iff\ $\tilde{\J} = \gf \J \gf^{-1}$. But this is precisely the T-duality relation for $(\J,\tilde\J)$.
\end{proof}

\end{document}